\documentclass[12pt]{article}
\usepackage{amsmath,amssymb,amsthm,amscd,graphicx}
\usepackage{latexsym}
\usepackage{enumerate}
\usepackage{amsbsy, amsfonts, textcomp}
\usepackage{makeidx}
\newtheorem{theorem}{Theorem}[]
\newtheorem{proposition}[theorem]{Proposition}
\newtheorem{lemma}[theorem]{Lemma}
\newtheorem{corollary}[theorem]{Corollary}

\theoremstyle{definition}

\newtheorem{example}[theorem]{Example}

\theoremstyle{remark}
\newtheorem{remark}[theorem]{Remark}
\newcommand{\G}{\mathbb{G}}

\newcommand{\C}{\mathbb{C}}

\newcommand{\N}{\mathbb{N}}

\def \rank{\operatorname{rank}}

\def \max{\operatorname{max}}

\begin{document}

\begin{center}

\Large

{\bf A new characterization of the invertibility \\ of polynomial maps}

\normalsize
\vspace{1cm}

EL\.ZBIETA ADAMUS \\ Faculty of Applied Mathematics, \\ AGH University of Science and Technology \\
al. Mickiewicza 30, 30-059 Krak\'ow, Poland \\
e-mail: esowa@agh.edu.pl \\

\vspace{0.5cm}
PAWE\L \ BOGDAN \\ Faculty of Mathematics and Computer Science, \\ Jagiellonian University \\
ul. \L ojasiewicza 6, 30-348 Krak\'ow, Poland \\
e-mail: pawel.bogdan@uj.edu.pl \\

\vspace{0.5cm}
TERESA CRESPO \\ Departament d'\`{A}lgebra i Geometria, \\ Universitat de Barcelona \\
Gran Via de les Corts Catalanes 585, 08007 Barcelona, Spain \\
e-mail: teresa.crespo@ub.edu

\vspace{0.5cm}
ZBIGNIEW HAJTO \\ Faculty of Mathematics and Computer Science, \\ Jagiellonian University \\
ul. \L ojasiewicza 6, 30-348 Krak\'ow, Poland \\
e-mail: zbigniew.hajto@uj.edu.pl

\end{center}

\vspace{0.2cm}

\begin{abstract}
In this paper we present an equivalent statement to the Jacobian conjecture. For a polynomial map $F$ on an affine space of dimension $n$,
we define recursively $n$ finite sequences of polynomials. We give an equivalent condition to the  invertibility of $F$ as well as a formula for
$F^{-1}$ in terms of these finite sequences of polynomials. Some examples illustrate the effective aspects of our approach. 
\end{abstract}

\section{Introduction}

The Jacobian Conjecture originated in the question raised by Keller in \cite{K} on the invertibility of polynomial maps with Jacobian determinant equal to~1. The question is still open in spite of the efforts of many mathematicians. We recall in the sequel the precise statement of the Jacobian Conjecture, some reduction theorems and other results we shall use. We refer to \cite{E} for a detailed account of the research on the Jacobian Conjecture and related topics.

Let $K$ be a field and $K[X]=K[X_1,\dots,X_n]$ the polynomial ring in the variables $X_1,\dots,X_n$ over $K$. A \emph{polynomial map} is a map $F=(F_1,\dots,F_n):K^n \rightarrow K^n$ of the form

$$(X_1,\dots,X_n)\mapsto (F_1(X_1,\dots,X_n),\dots,F_n(X_1,\dots,X_n)),$$

\noindent where $F_i \in K[X], 1 \leq i \leq n$. The polynomial map $F$ is \emph{invertible} if there exists a polynomial map $G=(G_1,\dots,G_n):K^n \rightarrow K^n$ such that $X_i=G_i(F_1,\dots,F_n),  1 \leq i \leq n$. We shall call $F$ a \emph{Keller map} if the Jacobian matrix

$$J=\left(\dfrac{\partial F_i}{\partial X_j}\right)_{\substack{1\leq i \leq n \\ 1\leq j \leq n}}$$

\noindent has determinant equal to 1. Clearly an invertible polynomial map $F$ has a Jacobian matrix $J$ with non zero determinant and may be transformed into a Keller map by composition with the linear automorphism with matrix $J(0)^{-1}$.

\vspace{0.5cm}
\noindent {\bf Jacobian Conjecture.} {\it Let $K$ be a field of characteristic zero. A Keller map $F:K^n \rightarrow K^n$ is invertible.}

\vspace{0.5cm}
In the sequel, $K$ will always denote a field of characteristic $0$.
For $F=(F_1,\dots,F_n) \in K[X]^n$, we define the \emph{degree} of $F$ as $\deg F= \max \{\deg F_i : 1\leq i \leq n\}$. It is known that if $F$ is a polynomial automorphism of $K^n$, then $\deg F^{-1} \leq (\deg F)^{n-1}$ (see \cite{BCW} or \cite{RW}).

The Jacobian conjecture for quadratic maps was proved by Wang in \cite{W}.  We state now the reduction of the Jacobian conjecture to the case of maps of third degree (see \cite{BCW}, \cite{Y}, \cite{D1} and \cite{D2}).

\begin{proposition}\label{red} \begin{enumerate}[a)] \item (Bass-Connell-Wright-Yagzhev) Given a Keller map $F:K^n \rightarrow K^n$, there exists a Keller map $\widetilde{F}:K^N \rightarrow K^N$, $N\geq n$ of the form $\widetilde{F}=Id+H$, where $H(X)$ is a cubic homogeneous map and having the following property: if $\widetilde{F}$ is invertible, then $F$ is invertible too.
\item (Dru\.zkowski) The cubic part $H$ may be chosen of the form

$$\left( (\sum_{j=1}^N a_{1j} X_j)^3, \dots,(\sum_{j=1}^N a_{Nj} X_j)^3 \right)$$

\noindent and with the matrix $A=(a_{ij})_{\substack{1\leq i \leq N \\ 1\leq j \leq N}}$ satisfying $A^2=0$.

\end{enumerate}

\end{proposition}

Polynomial maps in the Dru\.zkowski form are easier to handle than general cubic homogeneous polynomial maps. However we note the following result.

\begin{proposition}[\cite{E2} Proposition 2.9] Let $r \in \N$. If the Jacobian Conjecture holds for all cubic homogeneous polynomial maps in $r$ variables, then for all $n \in \N$  the Jacobian Conjecture holds for all polynomial maps of the form

$$F=X+(AX)^3$$

\noindent with $A \in \mathrm{M}_n(K)$ and $\rank A \leq r$.
\end{proposition}

In \cite{DR} Dru\.zkowski and Rusek give the following inversion formula for cubic homogeneous polynomial maps.

\begin{theorem}[\cite{DR}, Theorem 2.1] Let $H:K^n \rightarrow K^n$ be a cubic homogeneous polynomial map, $F=Id-H$ and let $G= \sum_{j=0}^{\infty} G_j$, where $G_j:K^n \rightarrow K^n$  is a homogeneous polynomial map of degree $j$, be the formal inverse of $F$. Then

$$\begin{array}{rll} G_1&=& Id, \\ [8pt] G_{2k+1} &=& \sum_{p+q+r=k-1} \varphi_H(G_{2p+1},G_{2q+1},G_{2r+1}), \, \forall k\geq 1, \\ [8pt] G_{2k} &=& 0 , \, \forall k\geq 1,
\end{array} $$

\noindent where $\varphi_H$ denotes the unique symmetric trilinear map such that $\varphi_H(X,X,X)=H(X)$.
\end{theorem}

\noindent As a corollary, they obtain that, if for some natural number $k$, we have

\begin{equation}\label{nest} G_{3^k +2}= \dots = G_{3^{k+1}}=0,
\end{equation}

\noindent then $F$ is a polynomial automorphism and $\deg F^{-1} \leq 3^k$. However, in \cite{GZ}, Gorni and Zampieri present an example of a
polynomial automorphism of $\C^4$ for which condition (\ref{nest}) is not satisfied for any $k$ (see example \ref{Aldo} below).


In this paper we present an algorithm providing a new characterization of the invertibility of
polynomial maps. Given a polynomial map $F:K^n \rightarrow K^n$ of
the form $F=Id+H$, where $H(X)$ has lower degree $\geq 2$, we define
recursively, for $1\leq i \leq n$, a sequence $P_k^i$ of
polynomials in $K[X]$ with $P_0^i=X_i$ such that $F$ is invertible
if and only if the alternating sum $\sum_{j=0}^{m-1}
(-1)^j P_j^i(X)$ satisfies a certain relation with $P_m^i$ for all
$i=1,\dots,n$, where $m$ is an integer given explicitly and depending
on the degrees of the components of $H$. When $F$ is invertible, its inverse
$F^{-1}$ is given in terms of these alternating sums of polynomials. In the
last section, we apply the algorithm to several examples of
polynomial maps, including the one of Gorni and Zampieri.

\section{A sufficient condition for invertibility}

Let us consider a polynomial map $F:K^n \rightarrow K^n$.
 Given a polynomial $P(X_1,\dots,X_n) \in K[X]=K[X_1,\dots,X_n]$, we define the following sequence of polynomials in $K[X]$,

$$\begin{array}{lll} P_0(X_1,\dots,X_n) &= &P(X_1,\dots,X_n), \\
P_1(X_1,\dots,X_n) &=& P_0(F_1,\dots,F_n)-P_0(X_1,\dots,X_n), \end{array}$$

\noindent and, assuming $P_{k-1}$ is defined,

$$P_k(X_1,\dots,X_n) = P_{k-1}(F_1,\dots,F_n)-P_{k-1}(X_1,\dots,X_n).$$

The following lemma is easy to prove.

\begin{lemma}\label{lem} For a positive integer $m$, we have

$$P(X_1,\dots,X_n)= \sum_{l=0}^{m-1} (-1)^l P_l(F_1,\dots,F_n)+(-1)^m P_m(X_1,\dots,X_n).$$

\noindent In particular, if we assume that for some integer $m$, $P_m(X_1,\dots,X_n)=0$, then

$$P(X_1,\dots,X_n)= \sum_{l=0}^{m-1} (-1)^l P_l(F_1,\dots,F_n).$$

\end{lemma}

\begin{corollary} Let $F:K^n \rightarrow K^n$ be a polynomial map. Let us consider the polynomial sequence $(P_k^i)$ constructed with $P=X_i$, $i=1,\dots, n$. Let us assume that for all $i=1,\dots, n$, there exists an integer $m_i$ such that $P_{m_i}^i = 0$. Then the inverse map $G$ of $F$ is given by

$$G_i(Y_1,Y_2,\dots,Y_n)= \sum_{l=0}^{m_i-1} (-1)^l P_l^i(Y_1,Y_2,\dots,Y_n), \, 1 \leq i \leq n.$$

\end{corollary}

The condition $P_{m_i}^i = 0$, for some integer $m_i$ for all $i=1,\dots, n$, is not necessary for the invertibility of $F$ (see example \ref{Aldo}). However we give in theorem \ref{gen} an equivalent condition to the invertibility of $F$ using a finite number of terms of the polynomial sequences $(P_k^i)$.
The following lemma gives a precise description of the polynomials $P_k^i$.

\begin{lemma}\label{lemgen}
  Let $F: K^n \rightarrow K^n$ be a polynomial map of the form
  \[\left\{ \begin{array}{l}
             F_1(X_1, \ldots, X_n) = X_1+H_1(X_1, \ldots, X_n)\\
              \qquad \qquad \vdots \\
             F_n(X_1, \ldots, X_n) = X_n+H_n(X_1, \ldots, X_n),
            \end{array} \right. \]
  where $H_i(X_1, \ldots, X_n)$ is a  polynomial in $X_1, \ldots, X_n$ of degree $D_i$ and lower degree $d_i$, with $ d_i\geq 2$, for $i=1, \ldots, n$. Let $d=\min d_i, D=\max D_i$.
Then for the polynomial sequence $(P^i_k)$ constructed with $P=X_i$ we have
that $P^i_k$ is a polynomial of degree $\leq D^{k-1}D_i$ and lower degree $\geq (k-1)(d-1)+d_i$.

In particular, if each $H_i$ is a homogeneous polynomial of degree $d$, we have

\begin{equation*}
P_k^i= \sum_{j=1}^{(d^k-1)/(d-1)-k+1} Q_{kj},
\end{equation*}

\noindent where $Q_{kj}$ is a homogeneous
polynomial in $X_1,\dots,X_n$ of degree $(k+j-1)(d-1)+1$.
\end{lemma}

\noindent {\it Proof.}  Let us consider, for a fixed $i$, the polynomial sequence

\begin{equation*}
\begin{array}{lll}
P_0^i(X_1,\dots,X_n) & = & X_i, \\
P_1^i(X_1,\dots,X_n) & = & F_i(X_1,\dots,X_n)-X_i=H_i(X_1,\dots,X_n), \\
P_2^i(X_1,\dots,X_n) & = & H_i(F_1,\dots,F_n)-H_i(X_1,\dots,X_n), \\
& \vdots &
\end{array}%
\end{equation*}

\noindent We write the Taylor series for the polynomial $H_i(F_1,\dots,F_n)=H_i(X_1+H_1,\dots,X_n+H_n)$ and obtain

\begin{equation*}
\begin{array}{lll}
P_2^i(X_1,\dots,X_n) & = & H_i(F_1,\dots,F_n)-H_i(X_1,\dots,X_n) \\
& = & Q_{21}^i+Q_{22}^i+\ldots+Q_{2D_i}^i%
\end{array}%
\end{equation*}

\noindent where

\[ Q_{21}^i  =  \sum_{j=1}^n  \frac{\partial H_i}{\partial X_j}H_j \]
\[ Q_{22}^i  =  \frac{1}{ 2!} \sum_{1\leq j_1,j_2 \leq n}^n \frac{ \partial^2 H_i }{ \partial X_{j_1} \partial X_{j_2}} H_{j_1} H_{j_2} \]
\vspace{5pt}
\[ \vdots\]
\[ Q^i_{ 2 D_i} = \frac{1}{D_i!} \sum_{j_1,j_2,\ldots, j_{D_i}=1}^n \frac{\partial ^{D_i}H_i}{\partial x_{j_1} \ldots \partial x_{j_{D_i}}}H_{j_1}\ldots H_{j_{D_i}}. \]

\noindent The polynomial $P_2^i$ has lower degree equal to the lower degree of $Q_{21}^i$, which is $\geq d+d_i-1$, and degree equal to the degree of $Q_{2D_i}^i$, which is  $\leq D \cdot D_i$. Let us prove by induction that
$P^i_k$ is a polynomial of degree $\leq D^{k-1}D_i$ and lower degree $\geq (k-1)(d-1)+d_i$. We have already seen it for $k=2$. Let us assume $P^i_{k-1}$ is a polynomial of degree $\leq D^{k-2}D_i$ and lower degree $\geq (k-2)(d-1)+d_i$.
We want to prove the property for $P_k^i$. We have

\begin{equation*}
P_k^i(X_1,\dots,X_n) = P_{k-1}^i(F_1,\dots,F_n)-P_{k-1}^i(X_1,\dots,X_n)
\end{equation*}

\noindent If $Q(X_1,\dots,X_n)$ is a polynomial of
degree $S$ and lower degree $s$,

\begin{equation*}
\begin{array}{l}
Q(F_1,\dots,F_n)-Q(X_1,\dots,X_n)=  \\ \quad \quad \quad \quad \quad \quad \sum_{j=1}^n \dfrac{\partial Q}{
\partial X_j} H_j + \ldots + \frac{1}{S!} \sum_{j_1, \ldots, j_S}^n \frac{\partial ^S Q}{\partial x_{j_1} \ldots, \partial x_{j_S}} H_{j_1} \ldots H_{j_S}
\end{array}
\end{equation*}

\noindent is a polynomial of degree $\leq S \cdot D$ and lower degree $\geq s-1+d$.
Hence $P_{k-1}^i(F_1,\dots,F_n)-P_{k-1}^i(X_1,\dots,X_n)$ is a polynomial
of degree $\leq D^{k-1}D_i$
and lower degree $\geq (k-1)(d-1)+d_i$.

The homogeneous case is proved analogously using induction. $\Box$

\begin{example}\label{Aldo} We shall consider the polynomial automorphism of $\C^4$ given in \cite{GZ} to prove that the condition $P_{m_i}^i = 0$, for some $m_i$, for all $i$, is not a necessary condition to the invertibility of $F$.

Let $p:=X_1X_3+X_2X_4$ and define $F$ by

$$\left\{ \begin{array}{lll} F_1 &=& X_1+pX_4 \\ F_2 &=& X_2-pX_3 \\F_3 &=& X_3+X_4^3 \\ F_4 &=& X_4 \end{array} \right.$$

\noindent Clearly $P_1^4=0$ and $P_2^3=0$. But $P_j^1$ and $P_j^2$
are not zero for any $j$. In order to prove that $P_j^1\neq 0$, we
shall prove by induction that the homogeneous summand of lowest
degree $Q_{j1}^1$ of $P_j^1$ has the following form depending on
the parity of $j$, for all $j \geq 2$.

$$\begin{array}{lll} Q_{2k,1}^1&=& X_1X_4^{4k} \\
Q_{2k+1,1}^1 &=&X_1X_3X_4^{4k+1} +X_2X_4^{4k+2} \end{array}$$

\noindent By calculation we obtain $Q_{21}^1=X_1X_4^4,
Q_{31}^1=X_1X_3X_4^5+X_2X_4^6$. Now, $Q_{2k,1}^1= X_1X_4^{4k}
\Rightarrow Q_{2k+1,1}^1=
X_4^{4k}H_1+4kX_1X_4^{4k-1}H_4=X_4^{4k}(X_1X_3+X_2X_4)X_4=X_1X_3X_4^{4k+1}
+X_2X_4^{4k+2}$ and $Q_{2k+1,1}^1 = X_1X_3X_4^{4k+1}
+X_2X_4^{4k+2} \Rightarrow Q_{2k+2,1}^1 =
X_3X_4^{4k+1}H_1+X_4^{4k+2}H_2+X_1X_4^{4k+1}
H_3+((4k+1)X_1X_3X_4^{4k}+(4k+2)X_2X_4^{4k+1})H_4$ \newline $=X_3X_4^{4k+1}(X_1X_3+X_2X_4)X_4-X_4^{4k+2}(X_1X_3+X_2X_4)X_3+X_1X_4^{4k+4}
=X_1X_4^{4(k+1)}.$

Analogously, in order to prove that $P_j^2\neq 0$, we shall prove
by induction that the homogeneous summand of lowest degree
$Q_{j1}^2$ of $P_j^2$ has the following form depending on the
parity of $j$, for all $j \geq 2$.

$$\begin{array}{lll} Q_{2k,1}^2&=& -2kX_1X_3X_4^{4k-1}-(2k-1)X_2X_4^{4k} \\
Q_{2k+1,1}^2&=&-X_1X_3^2X_4^{4k} -X_2X_3X_4^{4k+1}-2kX_1X_4^{4k+2}
\end{array}$$

\noindent By calculation we obtain $Q_{21}^2=-2X_1X_3X_4^3-X_2X_4^4,
Q_{31}^1=-X_1X_3^2X_4^4-X_2X_3X_4^5-2X_1X_4^6$. Now $Q_{2k,1}^2 =
-2kX_1X_3X_4^{4k-1}-(2k-1)X_2X_4^{4k} \Rightarrow Q_{2k+1,1}^2=
-2kX_3X_4^{4k-1}H_1-(2k-1)X_4^{4k}H_2-2kX_1X_4^{4k-1}H_3=-X_1X_3^2X_4^{4k}
-X_2X_3X_4^{4k+1}-2kX_1X_4^{4k+2}$ and
$Q_{2k+1,1}^2=-X_1X_3^2X_4^{4k} -X_2X_3X_4^{4k+1}-2kX_1X_4^{4k+2}$ \newline
$\Rightarrow Q_{2k+2,1}^2=(-X_3^2X_4^{4k}-2kX_4^{4k+2})H_1
-X_3X_4^{4k+1}H_2-(2X_1X_3X_4^{4k}+X_2X_4^{4k+1})H_3=-(2k+2)X_1X_3X_4^{4k+3}-(2k+1)X_2X_4^{4k+4}$.

\end{example}

\section{An equivalent condition to invertibility}

The following theorem gives an equivalent condition to the invertibility of $F$ using a finite number of terms in the polynomial sequences $(P_k^i)$.

 \begin{theorem}\label{gen}
    Let $F: K^n \rightarrow K^n$ be a polynomial map of the form
  \[\left\{ \begin{array}{lll}
             F_1(X_1, \ldots, X_n) &=& X_1+H_1(X_1, \ldots, X_n)\\
              & \vdots & \\
             F_n(X_1, \ldots, X_n) &=& X_n+H_n(X_1, \ldots, X_n),
            \end{array} \right. \]
  where $H_i(X_1, \ldots, X_n)$ is a  polynomial in $X_1, \ldots, X_n$ of degree $D_i$ and lower degree $d_i$, with $ d_i\geq 2$, for $i=1, \ldots, n$. Let $d=\min d_i, D=\max D_i$.
  The following conditions are equivalent:

\begin{enumerate}[1)]
\item $F$ is invertible.
\item For $i=1,\ldots, n$ and every $m > \frac{D^{n-1}-d_i}{d-1}+1$, we have
   \[ \sum_{j=0}^{m-1}(-1)^j P^i_j(X)=G_i(X)+R^i_m(X). \]
   where $G_i(X)$ is a polynomial of degree $ \leq D^{n-1}$, independent of $m$, and $R^i_m(X)$ is a polynomial satisfying $R^i_m(F)=(-1)^{m+1}P^i_m(X)$
   (with lower degree $\geq (m-1)(d-1)+d_i>D^{n-1}$).
\item For $i=1,\ldots, n$ and  $m = \lfloor \frac{D^{n-1}-d_i}{d-1}+1 \rfloor +1$, we have
   \[ \sum_{j=0}^{m-1}(-1)^j P^i_j(X)=G_i(X)+R^i_m(X). \]
   where $G_i(X)$ is a polynomial of degree $ \leq D^{n-1}$,
  and $R^i_m(X)$ is a polynomial satisfying $R^i_m(F)=(-1)^{m+1}P^i_m(X)$.
\end{enumerate}

  Moreover the inverse $G$ of $F$ is given by

  \[G_i(Y_1, \ldots, Y_n)=\sum_{l=0}^{m-1}(-1)^l\tilde{P}^i_l(Y_1, \ldots, Y_n), \, i=1, \ldots, n,\]
  where $\tilde{P}^i_l$ is the sum of homogeneous summands of $P^i_l$ of degree $ \leq D^{n-1}$ and $m$ is an integer $> \frac{D^{n-1}-d_i}{d-1}+1$.

 \end{theorem}

\noindent {\it Proof.} 1) $\Rightarrow$ 2): If $F$ is invertible,
then $G=F^{-1}$ has degree $\leq D^{n-1}$.
Applying lemma \ref{lem}, we obtain, for any positive integer $m$,

\begin{equation*}
X_i= \sum_{l=0}^{m-1} (-1)^l P_l(F_1,\dots,F_n)+(-1)^m P_m(X_1,\dots,X_n).
\end{equation*}

\noindent Since $X_i=G_i(F_1,\dots,F_n)$, we obtain the
following equality of polynomials in the variables $Y_1,\dots,Y_n$.

\begin{equation*}
\begin{array}{lll}
G_i(Y_1,\dots,Y_n) & = & \sum_{l=0}^{m-1} (-1)^l P_l(Y_1,\dots,Y_n) \\[8pt]
&  & +(-1)^m P_m(G_1(Y_1,\dots,Y_n),\dots,G_n(Y_1,\dots,Y_n)),
\end{array}
\end{equation*}

\noindent which implies

\begin{equation*}
\begin{array}{lll}
\sum_{l=0}^{m-1} (-1)^l P_l(Y_1,\dots,Y_n)& =  G_i(Y_1,\dots,Y_n)  \\[8pt]
& -(-1)^m P_m(G_1(Y_1,\dots,Y_n),\dots,G_n(Y_1,\dots,Y_n)),
\end{array}
\end{equation*}

\noindent Hence, writing

$$R_m^i(Y_1,\dots,Y_n):=-(-1)^m P_m(G_1(Y_1,\dots,Y_n),\dots,G_n(Y_1,\dots,Y_n)),$$

\noindent we obtain 2). Now, $G_i$ is a polynomial of degree at most $D^{n-1}$ in
$Y_1, \dots, Y_n$. For an integer $m$ such that $m > \frac{D^{n-1}-d_i}{d-1}+1$, $P_m$ is a
polynomial in the variables $X_1,\dots,X_n$  of lower degree bigger than $D^{n-1}$, hence the lower degree of $%
P_m(G_1(Y_1,\dots,Y_n),\dots,G_n(Y_1,\dots,Y_n))$ in the variables $%
Y_1,\dots,Y_n$ is bigger than $D^{n-1}$. Therefore, the sum of homogeneous
summands of degrees not bigger than $D^{n-1}$ in the righthand side of the
equality above is precisely $\sum_{l=0}^{m-1} (-1)^l \widetilde{P}
_l^i(Y_1,Y_2,\dots,Y_n)$.

\noindent 2) $\Rightarrow$ 3) is obvious.

\noindent 3) $\Rightarrow$ 1): Let us assume that for $m = \lfloor \frac{D^{n-1}-d_i}{d-1}+1 \rfloor +1$, we have

\begin{equation*} \sum_{j=0}^{m-1} (-1)^j P_j^i(X)=G_i(X) +R_m^i(X),
\end{equation*}

\noindent where $G_i(X)$ is a polynomial of degree $\leq D^{n-1}$ and $R_m^i(X)$ is a polynomial satisfying $R_m^i(F)=(-1)^{m+1} P_{m}^i(X)$. By lemma \ref{lem}, we have

\begin{equation*}
X_i= \sum_{l=0}^{m-1} (-1)^l P_l^i(F_1,\dots,F_n)+(-1)^m P_m^i(X_1,\dots,X_n).
\end{equation*}

\noindent We obtain then

\begin{equation*}\begin{array}{lll} X_i&=& G_i(F)+R_m^i(F)+(-1)^m P_m^i(X)\\ &=& G_i(F)+(-1)^{m+1} P_{m}^i(X)+(-1)^m P_m^i(X)\\ &=& G_i(F).\end{array}
\end{equation*}

\noindent Hence $F$ is invertible with inverse $G=(G_1,\dots,G_n)$. $\Box$

\vspace{0.5cm}

\section{Examples}

\subsection{} We consider the following nonhomogeneous Keller map in dimension 2.

$$\left\{ \begin{array}{lll} F_1 &= & X_1+(X_2+X_1^3)^2 \\ F_2 &=& X_2 + X_1^3 \end{array} \right.$$

\noindent Let us write $H_1:= (X_2+X_1^3)^2, H_2:=X_1^3$. With the notations in theorem \ref{gen}, we have $d_1=2, d_2=3, d=2, D=6$ and we obtain

$$\sum_{i=0}^5 (-1)^i P_i^1(X) = X_1-X_2^2+R_6^1(X),$$

\noindent where $R_6^1(X)$ is a polynomial of degree $6^5$ and lower degree 9 satisfying $R_6^1(F)=-P_6^1(X)$, and

$$\sum_{i=0}^4 (-1)^i P_i^2(X) = X_2-X_1^3+3 X_1^2 X_2^2-3 X_1 X_2^4+X_2^6+R_5^2(X),$$

\noindent where $R_5^2(X)$ is a polynomial of degree $3\cdot 6^3$ and lower degree 8 satisfying $R_5^2(F)=P_5^2(X)$. Hence the inverse of $F$ is given by

$$\left\{ \begin{array}{lll} G_1 &=& X_1-X_2^2\\ G_2 &=& X_2-X_1^3+3 X_1^2 X_2^2-3 X_1 X_2^4+X_2^6 \end{array} \right.$$

\subsection{} We consider the following Keller map $F$ in dimension 5.

$$\left\{ \begin{array}{lll} F_1 &= & X_1+a_1X_4^3+a_2x_4^2X_5+a_3X_4X_5^2+a_4x_5^3 +\dfrac{a_2c_5X_2X_4^2}{c_2}+\dfrac {2a_3c_5X_2X_4X_5}{c_2} \\[10pt] && +\dfrac{3a_4c_5X_2X_5^2}{c_2}-\dfrac{a_2e_2X_3X_4^2}{c_2}-\dfrac{2a_3e_2X_3X_4X_5}{c_2}-\dfrac{3a_4e_2X_3X_5^2}{c_2}\\[10pt] && +\dfrac {a_3c_5^2X_2^2X_4} {c_2^2}+\dfrac {3a_4c_5^2X_2^2X_5} {c_2^2}-\dfrac {2a_3c_5e_2X_2X_3X_4} {c_2^2}-\dfrac {6a_4c_5e_2X_2X_3X_5} {c_2^2} \\[10pt] && +\dfrac {a_3e_2^2X_3^2X_4} {c_2^2}+\dfrac {3a_4e_2^2X_3^3X_5} {c_2^2}+\dfrac {a_4e_5^3X_2^3}{c_2^3}-\dfrac {3a_4c_5^2e_2X_2^2X_3}{c_2^3}\\[10pt] &&+\dfrac {3a_4c_5e_2^2X_2X_3^2}{c_2^3}-\dfrac {a_4e_2^3X_3^3}{c_2^3} \\ [10pt] F_2 &=& X_2 + b_1X_4^3 \\ [10pt] F_3 &=& X_3+c_5X_2X_4^2+c_1X_4^3+c_2X_4^4X_5-e_2X_3X_4^2 \\ [10pt] F_4 &=& X_4 \\ [10pt] F_5 &=& X_5+e_2X_4^2X_5 +\dfrac {c_5e_2X_2X_4^2}{c_2} -\dfrac {e_2^2X_3X_4^2}{c_2}- \dfrac {(b_1c_5-c_1e_2)X_4^3}{c_2} \end{array} \right.$$

\noindent with parameters $a_1,a_2,a_3,a_4, b_1, c_1, c_2, c_5,
e_2$. By applying the algorithm we obtain $P_2^i=0$, for
$i=1,\dots,5$, hence $F$ is a quasi-translation, i.e.
$F^{-1}=2Id-F$.

\subsection{} We consider the following Keller map $F$ in dimension 6

$$\left\{ \begin{array}{lll}
F_1 &=& X_1+a_5e_1(X_1+X_2)^3/a_4+a_4X_2X_4X_6+a_5X_4X_5X_6 \\
F_2 &=& X_2-a_5e_1(X_1+X_2)^3/a_ 4 \\
F_3 &=& X_3+c_1X_1^3+c_2(X_1+X_5)^3+c_3(X_1+X_2)^3+c_4(X_1+X_4)^3+c_5X_6^3 \\
F_4 &=& X_4+d_4X_2X_6^2+a_5d_4X_5X_6^2/a_4 \\
F_5 &=& X_5+e_1(X_1+X_2)^3 \\
F_6 &=& X_6
\end{array} \right.$$

\noindent with parameters $a_4,a_5,c_1,c_2,c_3,c_4,c_5,d_4,e_1$. Denoting $G=F^{-1}$ and taking variables $(Y_1,\dots,Y_6)$ for $G$, we obtain $P_8^1=0$, $P_9^2=0$ and

\footnotesize

$$\begin{array}{l}G_1 = -(20a_5^4e_1 Y_6^9d_4^3Y_2^3Y_5^3a_4^3-6a_5^4e_1Y_6^4Y_4Y_5^3d_4Y_2a_4^2+a_5e_1a_4^6Y_6^9Y_2^6d_4^3+3Y_6^3a_5e_1Y_2^4d_4a_4^4\\  +3a_5e_1a_4^5Y_6^6Y_2^5d_4^2-a_5e_1a_4^6Y_2^3Y_4^3Y_6^3+a_4^5Y_2Y_4Y_6-Y_6^3a_4^5d_4Y_2^2+a_5^7e_1Y_6^9d_4^3Y_5^6 \\  +a_5e_1a_4^3Y_1^3+a_5e_1a_4^3Y_2^3-Y_1a_4^4-3a_5^2e_1a_4^5Y_2^2Y_4^3Y_6^3Y_5-6a_5e_1Y_1Y_2^2Y_4Y_6a_4^4 \\  +3a_5^5e_1Y_6^6Y_1d_4^2Y_5^4a_4+12a_5^4e_1Y_6^6Y_2^2d_4^2Y_5^3a_4^2-6a_5^4e_1Y_6^4Y_1Y_4d_4Y_5^3a_4^2+3a_5e_1a_4^3Y_1Y_2^2 \\ -
6a_5e_1a_4^5Y_6^4Y_1Y_2^3Y_4d_4+3Y_6^3a_5e_1Y_1^2Y_2^2d_4a_4^4+3a_5e_1a_4^5Y_1Y_2^2Y_4^2Y_6^2-30a_5^3e_1Y_6^7Y_4d_4^2Y_2^3Y_5^2a_4^4 \\  +6Y_6^3a_5e_1Y_1Y_2^3d_4a_4^4+
15a_5^3e_1Y_6^9d_4^3Y_2^4Y_5^2a_4^4+6a_5^2e_1a_4^5Y_6^9Y_2^5d_4^3Y_5-18a_5^2e_1Y_6^4Y_2^3Y_4d_4Y_5a_4^4\\  +12a_5^2e_1a_4^5Y_6^5Y_2^3Y_4^2d_4Y_5-
15a_5^2e_1a_4^5Y_6^7Y_2^4Y_4d_4^2Y_5+12a_5^2e_1Y_6^6Y_1Y_2^3d_4^2Y_5a_4^4-18a_5^2e_1Y_6^4Y_1Y_2^2Y_4d_4Y_5a_4^4 \\  -3a_5e_1a_4^6Y_6^7Y_2^5Y_4d_4^2+
3a_5e_1a_4^3Y_1^2Y_2+18a_5^3e_1Y_6^5Y_4^2d_4Y_2^2Y_5^2a_4^4+12a_5^2e_1Y_6^6Y_2^4d_4^2Y_5a_4^4\\  +3a_5e_1a_4^5Y_6^6Y_1Y_2^4d_4^2+3a_5e_1a_4^6Y_6^5Y_2^4Y_4^2d_4-
6a_5e_1Y_6^4a_4^5Y_2^4Y_4d_4+a_5Y_4Y_5Y_6a_4^4\\  +15a_5^5e_1Y_6^9d_4^3Y_2^2Y_5^4a_4^2+6a_5^6e_1Y_6^9Y_2d_4^3Y_5^5a_4-18a_5^3e_1Y_6^4Y_2^2Y_4d_4Y_5^2a_4^3+
18a_5^3e_1Y_6^6Y_1d_4^2Y_2^2Y_5^2a_4^3\\  -3a_5^6e_1Y_6^7Y_4Y_5^5d_4^2a_4-18a_5^3e_1Y_6^4Y_1Y_4d_4Y_2Y_5^2a_4^3+6a_5^3e_1Y_6^3Y_1d_4Y_5^2Y_2a_4^2-
30a_5^4e_1Y_6^7Y_4d_4^2Y_2^2Y_5^3a_4^3 \\  +12a_5^4e_1Y_6^5Y_2Y_4^2d_4Y_5^3a_4^3+12a_5^4e_1Y_6^6Y_1d_4^2Y_2Y_5^3a_4^2+3a_5^5e_1Y_6^6d_4^2Y_5^4Y_2a_4+
3a_5^3e_1Y_6^3Y_1^2d_4Y_5^2a_4^2 \\  +18a_5^3e_1Y_6^6Y_2^3d_4^2Y_5^2a_4^3+6a_5^2e_1Y_1Y_2Y_4^2Y_6^2Y_5a_4^4-3a_5e_1Y_2^3Y_4Y_6a_4^4-2Y_6^3Y_2a_5d_4Y_5a_4^4 \\  +
3a_5^5e_1Y_6^5Y_4^2Y_5^4d_4a_4^2-3a_5^3e_1Y_2Y_4^3Y_6^3Y_5^2a_4^4-3a_5e_1Y_1^2Y_2Y_4Y_6a_4^4+3a_5e_1a_4^5Y_2^3Y_4^2Y_6^2 \\  +6a_5^2e_1Y_2^2Y_4^2Y_6^2Y_5a_4^4-
3a_5^2e_1Y_4Y_5Y_6Y_2^2a_4^3-a_5^2Y_6^3d_4Y_5^2a_4^3-3a_5^2e_1Y_1^2Y_4Y_5Y_6a_4^3 \\  +3a_5^3e_1Y_4^2Y_5^2Y_6^2Y_2a_4^3+6Y_6^3a_5^2e_1Y_1^2Y_2d_4Y_5a_4^3-
15a_5^5e_1Y_6^7Y_2Y_4d_4^2Y_5^4a_4^2+12Y_6^3a_5^2e_1Y_1Y_2^2d_4Y_5a_4^3 \\  +6Y_6^3a_5^2e_1Y_2^3d_4Y_5a_4^3+3a_5^3e_1Y_6^3d_4Y_5^2Y_2^2a_4^2+
3a_5^3e_1Y_1Y_4^2Y_5^2Y_6^2a_4^3-6a_5^2e_1Y_1Y_4Y_5Y_6Y_2a_4^3 \\  -a_5^4e_1Y_4^3Y_5^3Y_6^3a_4^3)/a_4^4; \end{array}$$

$$\begin{array}{l}G_2=(20a_5^4e_1Y_6^9d_4^3Y_2^3Y_5^3a_4^3-6a_5^4e_1Y_6^4Y_4Y_5^3d_4Y_2a_4^2+a_5e_1a_4^6Y_6^9Y_2^6d_4^3+3Y_6^3a_5e_1Y_2^4d_4a_4^4\\
+3a_5e_1a_4^5Y_6^6Y_2^5d_4^2-
a_5e_1a_4^6Y_2^3Y_4^3Y_6^3+a_5^7e_1Y_6^9d_4^3Y_5^6+a_5e_1a_4^3Y_1^3\\ +a_5e_1a_4^3Y_2^3-3a_5^2e_1a_4^5Y_2^2Y_4^3Y_6^3Y_5-6a_5e_1Y_1Y_2^2Y_4Y_6a_4^4+
3a_5^5e_1Y_6^6Y_1d_4^2Y_5^4a_4\\ +12a_5^4e_1Y_6^6Y_2^2d_4^2Y_5^3a_4^2-6a_5^4e_1Y_6^4Y_1Y_4d_4Y_5^3a_4^2+3a_5e_1a_4^3Y_1Y_2^2-6a_5e_1a_4^5Y_6^4Y_1Y_2^3Y_4d_4 \\ +
3Y_6^3a_5e_1Y_1^2Y_2^2d_4a_4^4+3a_5e_1a_4^5Y_1Y_2^2Y_4^2Y_6^2-30a_5^3e_1Y_6^7Y_4d_4^2Y_2^3Y_5^2a_4^4+6Y_6^3a_5e_1Y_1Y_2^3d_4a_4^4 \\+
15a_5^3e_1Y_6^9d_4^3Y_2^4Y_5^2a_4^4+6a_5^2e_1a_4^5Y_6^9Y_2^5d_4^3Y_5-18a_5^2e_1Y_6^4Y_2^3Y_4d_4Y_5a_4^4+12a_5^2e_1a_4^5Y_6^5Y_2^3Y_4^2d_4Y_5 \\-
15a_5^2e_1a_4^5Y_6^7Y_2^4Y_4d_4^2Y_5+12a_5^2e_1Y_6^6Y_1Y_2^3d_4^2Y_5a_4^4-18a_5^2e_1Y_6^4Y_1Y_2^2Y_4d_4Y_5a_4^4-3a_5e_1a_4^6Y_6^7Y_2^5Y_4d_4^2 \\+
3a_5e_1a_4^3Y_1^2Y_2+18a_5^3e_1Y_6^5Y_4^2d_4Y_2^2Y_5^2a_4^4+12a_5^2e_1Y_6^6Y_2^4d_4^2Y_5a_4^4+3a_5e_1a_4^5Y_6^6Y_1Y_2^4d_4^2\\ +3a_5e_1a_4^6Y_6^5Y_2^4Y_4^2d_4-
6a_5e_1Y_6^4a_4^5Y_2^4Y_4d_4+15a_5^5e_1Y_6^9d_4^3Y_2^2Y_5^4a_4^2+6a_5^6e_1Y_6^9Y_2d_4^3Y_5^5a_4\\ -18a_5^3e_1Y_6^4Y_2^2Y_4d_4Y_5^2a_4^3+
18a_5^3e_1Y_6^6Y_1d_4^2Y_2^2Y_5^2a_4^3-3a_5^6e_1Y_6^7Y_4Y_5^5d_4^2a_4-18a_5^3e_1Y_6^4Y_1Y_4d_4Y_2Y_5^2a_4^3\\ +6a_5^3e_1Y_6^3Y_1d_4Y_5^2Y_2a_4^2-
30a_5^4e_1Y_6^7Y_4d_4^2Y_2^2Y_5^3a_4^3+12a_5^4e_1Y_6^5Y_2Y_4^2d_4Y_5^3a_4^3+12a_5^4e_1Y_6^6Y_1d_4^2Y_2Y_5^3a_4^2\\ +3a_5^5e_1Y_6^6d_4^2Y_5^4Y_2a_4+
3a_5^3e_1Y_6^3Y_1^2d_4Y_5^2a_4^2+18a_5^3e_1Y_6^6Y_2^3d_4^2Y_5^2a_4^3+6a_5^2e_1Y_1Y_2Y_4^2Y_6^2Y_5a_4^4\\ -3a_5e_1Y_2^3Y_4Y_6a_4^4+3a_5^5e_1Y_6^5Y_4^2Y_5^4d_4a_4^2-
3a_5^3e_1Y_2Y_4^3Y_6^3Y_5^2a_4^4-3a_5e_1Y_1^2Y_2Y_4Y_6a_4^4\\ +3a_5e_1a_4^5Y_2^3Y_4^2Y_6^2+6a_5^2e_1Y_2^2Y_4^2Y_6^2Y_5a_4^4-3a_5^2e_1Y_4Y_5Y_6Y_2^2a_4^3-
3a_5^2e_1Y_1^2Y_4Y_5Y_6a_4^3\\ +3a_5^3e_1Y_4^2Y_5^2Y_6^2Y_2a_4^3+6Y_6^3a_5^2e_1Y_1^2Y_2d_4Y_5a_4^3-15a_5^5e_1Y_6^7Y_2Y_4d_4^2Y_5^4a_4^2+
12Y_6^3a_5^2e_1Y_1Y_2^2d_4Y_5a_4^3\\ +6Y_6^3a_5^2e_1Y_2^3d_4Y_5a_4^3+3a_5^3e_1Y_6^3d_4Y_5^2Y_2^2a_4^2+3a_5^3e_1Y_1Y_4^2Y_5^2Y_6^2a_4^3-6a_5^2e_1Y_1Y_4Y_5Y_6Y_2a_4^3\\ -
a_5^4e_1Y_4^3Y_5^3Y_6^3a_4^3+Y_2a_4^4)/a_4^4.\end{array}$$

\normalsize

Now,

$$\begin{array}{lll} G_6 &=& Y_6 \\ G_5 &=& Y_5-e_1(G_1+G_2)^3 \\ G_4 &=& Y_4-d_4G_2G_6^2-a_5d_4G_5G_6^2/a_4 \\ G_3 &=& Y_3-c_1G_1^3-c_2 (G_1+G_5)^3-c_3(G_1+G_2)^3-c_4(G_1+G_4)^3+c_5G_6^3 .\end{array}$$

\subsection{} Let us consider again the polynomial automorphism of $\C^4$ given in example \ref{Aldo}. We have $p:=X_1X_3+X_2X_4$ and $F$ defined by

$$\left\{ \begin{array}{lll} F_1 &=& X_1+pX_4 \\ F_2 &=& X_2-pX_3 \\F_3 &=& X_3+X_4^3 \\ F_4 &=& X_4 \end{array} \right.$$

\noindent We obtain

$$\sum_{j=0}^{13} (-1)^j P_j^1(X)=X_1-X_1X_3X_4-X_2X_4^2+X_1X_4^4+R_{14}^1(X),$$

\noindent where

$$\begin{array}{lll}R_{14}^1(X)&=& -35X_1X_3X_4^{33}-10X_1X_3X_4^{45}-X_1X_3X_4^{49}-36X_1X_3X_4^{41}-6X_1X_3X_4^{29}\\&& -56X_1X_3X_4^{37}-X_1X_4^{28}-6X_2X_4^{30}
-15X_1X_4^{32}-35X_2X_4^{34}\\&&
-35X_1X_4^{36}-56X_2X_4^{38}-28X_1X_4^{40}-36X_2X_4^{42}-9X_1X_4^{44}\\
&& -10X_2X_4^{46}-X_1X_4^{48}-X_2X_4^{50}\end{array}$$

\noindent satisfies $R_{14}^1(F)+P_{14}^1(X)=0$. And

$$\sum_{j=0}^{13} (-1)^j P_j^2(X)= X_2-2X_1X_3X_4^3+X_2X_3X_4-X_2X_4^4+X_1X_4^6+X_1X_3^2+R_{14}^2(X),$$

\noindent where

$$\begin{array}{lll}R_{14}^2(X)&=&14X_1X_3X_4^{27}+6X_2X_3X_4^{29}+282X_1X_3X_4^{31}+6X_1X_3^2X_4^{28}+35X_2X_3X_4^{33}\\&& +910X_1X_3X_4^{35}
+35X_1X_3^2X_4^{32}+56X_2X_3X_4^{37}
+1064X_1X_3X_4^{39}+56X_1X_3^2X_4^{36}\\&&
+36X_2X_3X_4^{41}+558X_1X_3X_4^{43}+36X_1X_3^2X_4^{40}+12X_1X_3X_4^{51}+13X_2X_4^{28}\\&&
+77X_1X_4^{30}+267X_2X_4^{32}
+440X_1X_4^{34}+875X_2X_4^{36}+693X_1X_4^{38}\\ &&
+1036X_2X_4^{40}+440X_1X_4^{42}+549X_2X_4^{44}+121X_1X_4^{46}+133X_2X_4^{48}\\
&&  +12X_1X_4^{50}+12X_2X_4^{52} +10X_1X_3^2X_4^{44}+134
X_1X_3X_4^{47}+10X_2X_3X_4^{45}\\&&
+X_1X_3^2X_4^{48}+X_2X_3X_4^{49}\end{array}$$

\noindent satisfies $R_{14}^2(F)+P_{14}^2(X)=0$. Hence $G=F^{-1}$ is
given by

$$\left\{ \begin{array}{lll} G_1 &=& X_1-X_1X_3X_4-X_2X_4^2+X_1X_4^4 \\ G_2 &=& X_2-2X_1X_3X_4^3+X_2X_3X_4-X_2X_4^4+X_1X_4^6+X_1X_3^2 \\G_3 &=& X_3-X_4^3 \\ G_4 &=& X_4 \end{array} \right.$$

 \vspace{0.5cm} \noindent {\bf Acknowledgments.} E.
Adamus acknowledges support of the Polish Ministry of Science and
Higher Education. T. Crespo and Z. Hajto acknowledge support of
grant MTM2012-33830, Spanish Science Ministry.


\begin{thebibliography}{99}

\bibitem{BCW} H. Bass, E. Connell and D. Wright, \emph{The Jacobian Conjecture: Reduction of degree and formal expansion of the inverse}, Bull. Amer. Math. Soc. 7 (1982), 287-330.
\bibitem{D1} L. Dru\.zkowski, \emph{An effective approach to Keller's Jacobian conjecture}, Math. Ann. 264 (1983), 303-313.
\bibitem{D2} L. Dru\.zkowski, \emph{New reduction in the Jacobian conjecture}, Univ. Iagel. Acta Math.  No. 39  (2001), 203–-206.
\bibitem{DR} L. Dru\.zkowski, K. Rusek, \emph{The formal inverse and the Jacobian conjecture}, Ann. Polon. Math. 46 (1985), 85–-90.
\bibitem{E} A. van der Essen, \emph{Polynomial automorphisms and the Jacobian Conjecture}, Progress in Mathematics 190, Birkh\"auser Verlag, 2000.
\bibitem{E2} A. van der Essen, \emph{Seven lectures on polynomial automorphisms},  Automorphisms of affine spaces (Cura\c{c}ao, 1994), Kluwer Acad. Publ., Dordrecht, 1995, 3-39.
\bibitem{GZ} G. Gorni, G. Zampieri, \emph{Yagzhev polynomial mappings: on the structure of the Taylor expansion of their local inverse}, Ann. Polon. Math. 64 (1996), 285–-290.
\bibitem{K} O. H. Keller, \emph{Ganze Cremona-Transformationen}, Monatsh. Math. Phys. 47 (1939), 299-306.
\bibitem{RW} K. Rusek and T. Winiarski, \emph{Polynomial automorphisms of $\C^n$}, Univ. Iagell. Acta Math. 24 (1984), 143-149.
\bibitem{W} S. Wang, \emph{A Jacobian criterion for separability}, J. Algebra 65 (1980), 453-494.
\bibitem{Y} A. V. Yagzhev, \emph{On a problem of O.H. Keller}, Sibirsk. Mat. Zh. 21 (1980), 747-754.
\end{thebibliography}
\end{document}